\RequirePackage{etoolbox}
\csdef{input@path}{{style/}{graphics/}}
\documentclass[MSNbibl,number,citesort,seceqn,dvips]{arxbj}
\usepackage{upgreek}

\aid{0}
\volume{21}
\issue{1}
\pubyear{2015}
\firstpage{604}
\lastpage{617}
\doi{10.3150/13-BEJ581} 

\makeatletter

\newtheorem{theorem}{Theorem}[section]
\newtheorem{lemma}[theorem]{Lemma}
\newtheorem{corollary}[theorem]{Corollary}
\newremark{remark}[theorem]{Remark}
\newremark{quasi}{Quasi-infinitely divisible distribution}
\newremark{wkquasi}{Pretended-infinitely divisible distribution}

\newcommand{\eqref}[1]{(\ref{#1})}
\renewcommand{\emptyset}{\varnothing}
\newcommand{\ID}{\operatorname{ID}}
\makeatother

\begin{document}
\begin{frontmatter}

\title{A complete Riemann zeta distribution and the Riemann hypothesis}
\runtitle{A complete Riemann zeta distribution and the Riemann hypothesis}

\begin{aug}
\author{\inits{T.}\fnms{Takashi} \snm{Nakamura}\ead
[label=e1]{nakamura\_takashi@ma.noda.tus.ac.jp}}
\address[]{Department of Mathematics Faculty of Science and
Technology, Tokyo University of Science, Noda, Chiba 278-8510, Japan.
\printead{e1}}
\end{aug}

\received{\smonth{7} \syear{2013}}
\revised{\smonth{11} \syear{2013}}

%
\begin{abstract}
Let $\sigma, t \in{\mathbb{R}}$, $s=\sigma+\mathrm{{i}}t$, $\Gamma
(s)$ be the Gamma function,
$\zeta(s)$ be the Riemann zeta function and $\xi(s) := s(s-1)\pi
^{-s/2} \Gamma(s/2) \zeta(s)$
be the complete Riemann zeta function. We show that $\Xi_\sigma(t):=
\xi
(\sigma- \mathrm{{i}}t)/\xi(\sigma)$
is a characteristic function for any $\sigma\in{\mathbb{R}}$ by
giving the probability density function.
Next we prove that the Riemann hypothesis is true if and only if each
$\Xi_\sigma(t)$ is a pretended-infinitely
divisible characteristic function, which is defined in this paper, for
each $1/2 < \sigma<1$. Moreover,
we show that $\Xi_\sigma(t)$ is a pretended-infinitely divisible
characteristic function when $\sigma=1$.
Finally we prove that the characteristic function $\Xi_\sigma(t)$ is
not infinitely divisible but
quasi-infinitely divisible for any $\sigma>1$.
\end{abstract}

%
\begin{keyword}
\kwd{characteristic function}
\kwd{L\'evy--Khintchine representation}
\kwd{Riemann hypothesis}
\kwd{zeta distribution}
\end{keyword}

\end{frontmatter}
%
\section{Introduction and main results}
\subsection{Riemann zeta function and distribution}\label{sec1.1}
The famous Riemann zeta function $\zeta(s)$ is a function of a complex
variable $s=\sigma+\mathrm{{i}}t$, for $\sigma>1$ defined by
\begin{eqnarray*}
\zeta(s) := \sum_{n=1}^{\infty}
\frac{1}{n^s} = \prod_p \biggl( 1 -
\frac
{1}{p^s} \biggr)^{-1},
\end{eqnarray*}
where the letter $p$ is a prime number, and the product of $\prod_p$ is
taken over all primes. The Dirichlet series $\sum_{n=1}^{\infty}
n^{-s}$ and the Euler product $\prod_p ( 1 - p^{-s} )^{-1}$ converges
absolutely in the half-plane $\sigma>1$ and uniformly in each compact
subset of this half-plane. The Riemann zeta function is a meromorphic
function on the whole complex plane, which is holomorphic everywhere
except for a simple pole at $s = 1$ with residue $1$. Denote the Gamma
function by $\Gamma(s)$. We have the following functional equation of
the complete Riemann zeta function $\xi(s)$ (see, for example,
Titchmarsh \cite{Tit}, (2.1.13))
%
\begin{equation}
\label{eq:feq} \xi(s) = \xi(1-s), \qquad \xi(s) :=s(s-1) \uppi^{-s/2} \Gamma
\biggl( \frac{s}{2} \biggr) \zeta(s) .
\end{equation}

In view of the Euler product, it is seen easily that $\zeta(s)$ has no
zeros in the half-plane $\sigma>1$. It follows from the functional
equation (\ref{eq:feq}) and basic properties of the Gamma-function that
$\zeta(s)$ vanishes in $\sigma< 0$ exactly at the so-called trivial
zeros $s=-2m$, $m \in{\mathbb{N}}$. In 1859, Riemann stated that it
seems likely that all nontrivial zeros lie on the so-called critical
line $\sigma= 1/2$. This is the famous, yet unproved Riemann
hypothesis. In 1896, Hadamard and de la Vall\'ee-Poussin independently proved that $\zeta(1+\mathrm{{i}}t) \ne0$ for any $t
\in{\mathbb{R}}$ (see Titchmarsh \cite{Tit}, page 45). Hence, we can also see that no zeros of $\zeta(s)$
lie on the line $\Re(s) = 0$ by (\ref{eq:feq}). Therefore, the Riemann
hypothesis is rewritten equivalently as
\[
\mbox{\textit{Riemann hypothesis}} \quad \zeta(s) \ne0 \qquad \mbox{for } 1/2 < \sigma<
1.
\]

Put $Z_{\sigma}(t):=\zeta(\sigma-\mathrm{i}t)/\zeta(\sigma)$,
$t\in
{\mathbb{R}}$, then $Z_{\sigma}(t)$ is known to be a characteristic
function when $\sigma>1$ (see Khintchine \cite{Khi} or
Gnedenko and Kolmogorov \cite{GK68}, page 75). A
distribution $\mu_{\sigma}$ on ${\mathbb{R}}$ is said to be a Riemann
zeta distribution with parameter $\sigma$ if it has $Z_{\sigma}(t)$ as
its characteristic function. Recently, the Riemann zeta distribution is
investigated by Lin and Hu \cite{Lin}, and Gut
\cite{Gut}. On the other hand, in Aoyama and Nakamura \cite{ANpe}, Remark~1.13, it is showed that $Z_{\sigma}(t)$ is not a
characteristic function for any $1/2 \le\sigma\le1$. Afterwards,
Nakamura \cite{Nakamr} showed that $F_\sigma(t)$, where
$F_\sigma(t):= f_\sigma
(t)/f_\sigma(0)$ and $f_\sigma(t):=\zeta(\sigma-\mathrm
{{i}}t)/(\sigma
-\mathrm{{i}}t)$, is a characteristic function for any $0 < \sigma\ne1$.

Note that there are some other papers connected to Riemann zeta
function in probabilistic view. Biane Pitman and Yor \cite{BPY} reviewed
known results about
$\xi(s)$ which are related to one-dimensional Brownian motion and to
higher dimensional Bessel processes. Lagarias and Rains \cite{LR}
treated $\uppi^{-s/2}
\Gamma(s/2) \zeta(s)$ and its generalizations and gave results
connected to infinite divisibility.

\subsection{Infinitely divisible and quasi-infinitely divisible distributions}
A probability measure $\mu$ on ${\mathbb{R}}$ is infinitely divisible
if, for any positive integer $n$, there is a probability measure $\mu
_n$ on ${\mathbb{R}}$ such that $\mu=\mu_n^{n*}$, where $\mu
_n^{n*}$ is
the $n$-fold convolution of $\mu_n$. For instance, normal, degenerate,
Poisson and compound Poisson distributions are infinitely divisible.

Let $\widehat{\mu}(t)$ be the characteristic function of a probability
measure $\mu$ on ${\mathbb{R}}$ and $\ID({\mathbb{R}})$ be the class of
all infinitely divisible distributions on ${\mathbb{R}}$. The following
L\'evy--Khintchine representation is well known (see Sato
\cite{S99}, Section~2). Put $D_b:=\{x \in{\mathbb{R}}\dvt  -b \le x \le
b\}$, where
$b>0$. If $\mu\in \ID({\mathbb{R}})$, then one has
%
\begin{equation}
\label{INF} \widehat{\mu}(t) = \exp \biggl[-\frac{a}{2} t^2 +
\mathrm{i} \lambda t +\int_{\mathbb{R}} \bigl(\mathrm{e}^{\mathrm{i} tx}-1-
\mathrm{i} tx 1_{D_b}(x) \bigr)\nu(\mathrm{d}x) \biggr],\qquad  t\in{\mathbb{R}},
\end{equation}
where $a \ge0$, $\lambda\in{\mathbb{R}}$ and $\nu$ is a measure on
${\mathbb{R}}$ satisfies $\nu(\{0\}) = 0$ and $\int_{\mathbb{R}}
(|x|^{2} \wedge1) \nu(\mathrm{d}x) < \infty$. Moreover, the representation of
$\widehat{\mu}$ in (\ref{INF}) by $a, \nu$, and $\lambda$ is
unique. If
the L\'evy measure $\nu$ in $\eqref{INF}$ satisfies $\int_{|x|<1}|x|\nu
(\mathrm{d}x)<\infty$, then $\eqref{INF}$ can be written by
%
\begin{equation}
\label{INF2} \widehat{\mu}(t) = \exp \biggl[-\frac{a}{2} t^2 +
\mathrm{i} \lambda _0 t + \int_{\mathbb{R}}
\bigl(\mathrm{e}^{\mathrm{i} tx}-1 \bigr)\nu(\mathrm{d}x) \biggr], \qquad \lambda _0 \in{
\mathbb{R}}.
\end{equation}

For example, the L\'evy measure of $Z_{\sigma}(t):=\zeta(\sigma
-\mathrm{i}t)/\zeta(\sigma)$ can be given as in the following (see
Gnedenko and Kolmogorov \cite{GK68}, page 75). Let $\delta_x$ be the delta
measure at $x$. Then
we have
%
\begin{equation}
\label{eq:rzlm1} \log Z_{\sigma}(t) = \int_0^{\infty}
\bigl(\mathrm{e}^{\mathrm{i}tx}-1 \bigr)N_{\sigma}(\mathrm{d}x),\qquad  N_{\sigma}(\mathrm{d}x):=\sum
_{p}\sum_{r=1}^{\infty}
\frac{p^{-r\sigma
}}{r}\delta _{r\log p}(\mathrm{d}x).
\end{equation}

On the other hand, there are non-infinitely divisible distributions
whose characteristic functions are the quotients of two infinitely
divisible characteristic functions.
That class is called class of \textit{quasi-infinitely divisible
distributions} and is defined as follows.
%
\begin{quasi*}
A distribution $\mu$ on ${\mathbb{R}}$ is called \textit
{quasi-infinitely divisible} if it has a form of \eqref{INF} with $a
\in{\mathbb{R}}$ and the corresponding measure $\nu$ is a signed
measure on ${\mathbb{R}}$ with total variation measure $|\nu|$
satisfying $\nu(\{0\}) =0$ and $\int_{\mathbb{R}} (|x|^{2} \wedge1)
|\nu|(\mathrm{d}x) < \infty$.
\end{quasi*}

We have to mention that the triplet $(a,\nu,\lambda)$ in this case is
also unique if each component exists and that infinitely divisible
distributions on ${\mathbb{R}}$ are quasi-infinitely divisible if and
only if $a \ge0$ and the negative part of $\nu$ in the Jordan
decomposition equals zero.
The measure $\nu$ is called \textit{quasi}-L\'evy measure and has
appeared in some books and papers, for example, Gnedenko and Kolmogorov
\cite{GK68}, page 81,
Lindner and Sato \cite{LSm}, Niedbalska-Rajba \cite{N-R}, and
others (see also Sato \cite{Sato12}, Section~2.4).

\subsection{Main results}\label{sec1.3}
In the present paper, we give a complete Riemann zeta distribution by
the normalized complete Riemann zeta function
\begin{eqnarray*}
\Xi_\sigma(t) &:= &\frac{\xi(\sigma- \mathrm{{i}}t)}{\xi(\sigma
)},
\\
\xi(\sigma- \mathrm{{i}}t)&:=& (\sigma- \mathrm{{i}}t) (\sigma-1 - \mathrm{{i}}t)
\uppi^{(\mathrm{{i}}t-\sigma)/2} \Gamma \biggl(\frac
{\sigma
-\mathrm{{i}}t}{2} \biggr) \zeta(\sigma-
\mathrm{{i}}t),
\end{eqnarray*}
for any $\sigma\in{\mathbb{R}}$. It should be mentioned that $\Xi
_\sigma(t)$ is symmetric about the vertical axis $\sigma= 1/2$ by the
functional equation (\ref{eq:feq}). Therefore, we only have to consider
the case $\sigma\ge1/2$. In order to state the main results, we
introduce the following pretended-infinitely divisible distribution.
%
\begin{wkquasi*}
A distribution $\mu$ on ${\mathbb{R}}$ is called \textit
{pretended-infinitely divisible} if it has a form of \eqref{INF} with
$a \in{\mathbb{R}}$ and the corresponding measure $\nu$ is a signed
measure on ${\mathbb{R}}$ with $\nu(\{0\}) =0$.
\end{wkquasi*}

Namely, pretended-infinitely divisible distributions are infinitely
divisible or quasi-infinitely divisible distributions without the
condition $\int_{\mathbb{R}} (|x|^{2} \wedge1) |\nu|(\mathrm{d}x) < \infty$.

The main results in this paper are following four theorems.

\begin{theorem}\label{th:1}
The function $\Xi_\sigma(t)$ is a characteristic function for any
$\sigma\in{\mathbb{R}}$. Moreover, the probability density function
$P_\sigma(y)$ is given as follows:
%
\begin{eqnarray}
\label{eq:defP} P_\sigma(y) := %
\cases{  \displaystyle\frac{2}{\xi(\sigma)}
\displaystyle \sum_{n=1}^\infty f\bigl(n\mathrm{e}^{-y}
\bigr) \mathrm{e}^{-\sigma y}, &\quad  $y \le0$,
\cr
 \displaystyle\frac{2}{\xi(\sigma)} \displaystyle \sum
_{n=1}^\infty f\bigl(n\mathrm{e}^{y}\bigr)
\mathrm{e}^{(1-\sigma)y}, &\quad  $y>0$, } %
\end{eqnarray}
where $f (x) := 2\uppi(2\uppi x^4-3x^2)\mathrm{e}^{-\uppi x^2}$.
\end{theorem}

Let ${\mathcal{Z}}$ and ${\mathcal{Z}}_+$ be the set of zeros of the
Riemann zeta function which lie in the critical strip $\{s\in{\mathbb
{C}} \dvt  0 < \Re(s) <1\}$, and the region $\{s\in{\mathbb{C}} \dvt  0 <
\Re
(s) <1, \Im(s) >0 \}$, respectively. If the Riemann hypothesis is
true, then each $\rho\in{\mathcal{Z}}_+$ can be expressed by $\rho=
1/2 + \mathrm{{i}}\gamma$, where $\gamma>0$.

\begin{theorem}\label{th:2}
The characteristic function $\Xi_\sigma(t)$ is a pretended-infinitely
divisible characteristic function for any $1/2 < \sigma<1$ if and only
if the Riemann hypothesis is true. Furthermore, we have
%
\begin{eqnarray}
\label{eq:th2} \Xi_\sigma(t) &=& \exp \biggl[ \int_0^\infty
\bigl(\mathrm{e}^{\mathrm{{i}}tx}-1 \bigr)\nu _{\sigma
}(\mathrm{d}x) \biggr],\nonumber\\[-8pt]\\[-8pt]
\nu_{\sigma}(\mathrm{d}x)&:=& -\sum_{1/2+\mathrm{{i}}\gamma\in{\mathcal{Z}}_+}
\frac{2\cos
(\gamma
x)}{x\mathrm{e}^{(\sigma-1/2)x}} (\mathrm{d}x) ,\nonumber
\end{eqnarray}
under the Riemann hypothesis.
\end{theorem}

Let ${\mathcal{Z}}_+^R$ be the set of zeros of $\zeta(s)$ which lie on
the half line $\{s\in{\mathbb{C}} \dvt \Re(s) =1/2, \Im(s) >0 \}$ and
${\mathcal{Z}}_+^N$ be the set of zeros of $\zeta(s)$ which lie in the
region $\{s\in{\mathbb{C}} \dvt  1/2 < \Re(s) <1, \Im(s) >0 \}$. Note
that ${\mathcal{Z}}_+^N = \emptyset$ if and only if the Riemann
hypothesis is true. One has ${\mathcal{Z}}= \{ \rho, 1-\rho\dvt  \rho
\in
{\mathcal{Z}}_+^R \} \cup\{ \rho, 1-\rho, \overline{\rho},
1-\overline
{\rho}\dvt  \rho\in{\mathcal{Z}}_+^N \}$ from $\xi(s)= \xi(1-s)$ and
$\xi(\overline{s})= \overline{\xi(s)}$.

\begin{theorem}\label{th:3}
When $\sigma\ge1$, we have
%
\begin{eqnarray}
\label{eq:th3} \Xi_\sigma(t)& =& \exp \biggl[ \int_0^\infty
\bigl(\mathrm{e}^{\mathrm{{i}}tx}-1 \bigr)\nu_{\sigma}(\mathrm{d}x) \biggr],\nonumber
\\[-8pt]\\[-8pt]
\nu_{\sigma}(\mathrm{d}x)&:=& - \sum_{1/2+\mathrm{{i}}\gamma\in{\mathcal{Z}}_+^R}
\frac{2\cos(\gamma x)}{x\mathrm{e}^{(\sigma-1/2)x}} (\mathrm{d}x) - \sum_{\beta+\mathrm{{i}}\gamma\in{\mathcal{Z}}_+^N} \biggl(
\frac{2\cos(\gamma x)}{x\mathrm{e}^{(\sigma-\beta)x}} + \frac{2\cos(\gamma
x)}{x\mathrm{e}^{(\sigma-1+\beta)x}} \biggr) (\mathrm{d}x).\nonumber
\end{eqnarray}
Especially, $\Xi_\sigma(t)$ is a pretended-infinitely divisible
characteristic function when $\sigma=1$.
\end{theorem}

\begin{theorem}\label{th:4}
When $\sigma>1$, we have
\begin{eqnarray*}
\Xi_\sigma(t)& =& \exp \biggl[ \mathrm{i}t \lambda_\sigma+ \int
_0^\infty \bigl(\mathrm{e}^{\mathrm{{i}}tx}-1-\mathrm
{{i}}tx1_{D_{1/2}}(x) \bigr)\nu_{\sigma}(\mathrm{d}x) \biggr],
\\
\lambda_\sigma&:= &\frac{\mathrm{e}^{-\sigma/2}-1}{\sigma} + \frac
{\mathrm{e}^{(1-\sigma
)/2}-1}{\sigma-1}+
\frac{\log\uppi}{2}
 \\
&&{}+ \frac{1}{2} \int_0^1
\biggl( \frac{\mathrm{e}^{-\sigma x/2}}{1-\mathrm{e}^{-x}} - \frac
{\mathrm{e}^{-x}}{x} \biggr)\, \mathrm{d}x - \frac{1}{2}
\int_1^\infty \mathrm{e}^{-x} \frac{\mathrm{d}x}{x},
\\
\nu_{\sigma}(\mathrm{d}x)&:=& \frac{1(\mathrm{d}x)}{x \mathrm{e}^{\sigma x}(1-\mathrm{e}^{-2x})}-\frac
{1+\mathrm{e}^x}{x\mathrm{e}^{\sigma x}} (\mathrm{d}x) + \sum
_{p}\sum_{r=1}^{\infty}
\frac{p^{-r\sigma}}{r}\delta_{r\log
p}(\mathrm{d}x) .
\end{eqnarray*}
Therefore, the characteristic function $\Xi_\sigma(t)$ is not
infinitely divisible but quasi-infinitely divisible when $\sigma>1$.
\end{theorem}

We call the distribution defined by the characteristic function $\Xi
_\sigma(t)$ the \textit{{completed Riemann zeta distribution}}. It is
well known that $\zeta(s)$ has zeros on $\Re(s) =1/2$ (see
Titchmarsh \cite{Tit}, Section~10). By the definition of
pretended-infinitely divisible
distribution and the fact that $\exp(z) \ne0$ for any $z \in{\mathbb
{C}}$, the characteristic function does not have zeros. Thus, $\Xi
_\sigma(t)$ is not even a pretended-infinitely divisible
characteristic function when $\sigma=1/2$.

\section{Proofs}
\subsection{Proof of Theorem \texorpdfstring{\protect\ref{th:1}}{1.1}}
We quote the following fact from Patterson \cite{Pat} (see
also Biane Pitman and Yor \cite{BPY}, Section 2).

\begin{lemma}[(see {Patterson \cite{Pat}}, Section 2.10)]
Let $f (x) := 2\uppi(2\uppi x^4-3x^2)\mathrm{e}^{-\uppi x^2}$. Then we have
%
\begin{eqnarray}
\label{eq:tit1} \xi(s) = 2 \int_1^\infty\sum
_{n=1}^\infty f(nx) \bigl( x^{s-1/2} +
x^{1/2-s} \bigr) x^{-1/2}\,\mathrm{d}x .
\end{eqnarray}
Note that the last integral is absolutely convergent for all values of $s$.
\end{lemma}

\begin{pf*}{Proof of Theorem~\ref{th:1}}
By (\ref{eq:tit1}) and the change of variables $x=\mathrm{e}^{-y}$ and
$x=\mathrm{e}^{y}$, we have
\begin{eqnarray*}
\xi(\sigma- \mathrm{{i}}t) &=& 2 \int_1^\infty\sum
_{n=1}^\infty f(nx) x^{\sigma- \mathrm{{i}}t-1} \,\mathrm{d}x + 2
\int_1^\infty\sum_{n=1}^\infty
f(nx) x^{\mathrm{{i}}t-\sigma} \,\mathrm{d}x
\\
&= & 2 \int_0^{-\infty} \sum
_{n=1}^\infty f\bigl(n\mathrm{e}^{-y}\bigr)
\mathrm{e}^{(1+\mathrm{{i}}t
-\sigma)y} \bigl(-\mathrm{e}^{-y}\bigr)\,\mathrm{d}y + 2 \int_0^\infty
\sum_{n=1}^\infty f\bigl(n\mathrm{e}^y\bigr)
\mathrm{e}^{(\mathrm{{i}}t-\sigma)y} \bigl(\mathrm{e}^{y}\bigr)\,\mathrm{d}y
\\
&= & 2 \int_{-\infty}^0 \mathrm{e}^{\mathrm{{i}}ty} \sum
_{n=1}^\infty f\bigl(n\mathrm{e}^{-y}\bigr)
\mathrm{e}^{-\sigma y} \,\mathrm{d}y + 2 \int_0^\infty
\mathrm{e}^{\mathrm{{i}}ty} \sum_{n=1}^\infty f
\bigl(n\mathrm{e}^{y}\bigr) \mathrm{e}^{(1-\sigma)y}\, \mathrm{d}y.
\end{eqnarray*}
Obviously, we have $f (x) = 2\uppi(2\uppi x^4-3x^2)\mathrm{e}^{-\uppi x^2} > 0$ for
any $x\ge1$. Hence, one has $f(n\mathrm{e}^{-y})>0$ for any $y \le0$ and $n
\in{\mathbb{N}}$, and $f(n\mathrm{e}^y)>0$ for any $y > 0$ and $n \in{\mathbb
{N}}$. Thus it holds that
\[
\sum_{n=1}^\infty f\bigl(n\mathrm{e}^{-y}
\bigr) \mathrm{e}^{-\sigma y} >0,\qquad  y\le0 \quad \mbox{and}\quad  \sum_{n=1}^\infty
f\bigl(n\mathrm{e}^{y}\bigr) \mathrm{e}^{(1-\sigma)y}>0,\qquad  y> 0.
\]
On the other hand, we have
\[
\xi(\sigma) = 2 \int_{-\infty}^0 \sum
_{n=1}^\infty f\bigl(n  \mathrm{e}^{-y}\bigr)
\mathrm{e}^{-\sigma y} \,\mathrm{d}y + 2 \int_0^\infty\sum
_{n=1}^\infty f\bigl(n\mathrm{e}^{y}\bigr)
\mathrm{e}^{(1-\sigma)y}\, \mathrm{d}y >0
\]
from (\ref{eq:tit1}) and the argument above. Hence, $P_\sigma(y)$
defined by (\ref{eq:defP}) is nonnegative. Therefore, we have $\Xi
_\sigma(t)= \int_{\mathbb{R}} \mathrm{e}^{\mathrm{{i}}ty}P_\sigma(y)\,\mathrm{d}y$, where
$P_\sigma(y)$ is the probability density function.
\end{pf*}

\begin{remark}
It should be emphasised that $\Xi_\sigma(t)$ is a characteristic
function for any $\sigma\in{\mathbb{R}}$. On the other hand,
$F_\sigma
(t):= f_\sigma(t)/f_\sigma(0)$, where $f_\sigma(t):=\zeta(\sigma
-\mathrm
{{i}}t)/(\sigma-\mathrm{{i}}t)$, is not a characteristic function for
$\sigma= 0, 1$ and $\sigma< -1/2$. This is proved as follows. When
$\sigma=1$, it is well known that $\zeta(1+\mathrm{{i}}t) \ne0$, $t
\ne0$, and $\zeta(s)$ has an only one pole at $s=1$. Hence, we have
\[
F_1 (t) = \frac{1}{\zeta(1)} \frac{\zeta(1+\mathrm
{{i}}t)}{1+\mathrm
{{i}}t} = 0 \qquad \mbox{for any } t
\ne0,
\]
which contradicts the uniform continuity of characteristic function
$\widehat{\mu}(t)$ and $\widehat{\mu}(0)=1$. A~similar argument can be
done when $\sigma=0$ since $\zeta(s)/ s$ has a simple pole at $s=0$.
By (\ref{eq:feq}) and Stirling's formula, one has
\[
\bigl |\zeta(s)\bigr | = \uppi^{\sigma-1/2} \bigl(|t/2|+2\bigr)^{-\sigma+1/2} \bigl( 1 + \mathrm{O}
\bigl(\bigl(|t|+2\bigr)^{-1}\bigr) \bigr) \bigl |\zeta(1-s)\bigr |
\]
for $\sigma<0$. On the other hand, for any $\varepsilon>0$ there are
arbitrarily large $t$ which satisfy $|\zeta(\sigma+ \mathrm{{i}}t)| >
(1- \varepsilon) \zeta(\sigma)$ when $\sigma>1$ (see Titchmarsh \cite{Tit}, Theorem~8.4). Thus, we can find $t$ which satisfies $|\zeta(s)| > \uppi
^{\sigma-1/2} |t/2|^{-\sigma+1/2} \zeta(1-\sigma)/2$. Hence, there
exists $t \in{\mathbb{R}}$ such that $|F_\sigma(t)| >1$ when $\sigma<
-1/2$ by the factor $|t/2|^{-\sigma+1/2}$.
\end{remark}

The absolute value of a characteristic function is not greater than $1$
(see for instance Sato \cite{S99}, Proposition~2.5).
Hence, we have the
following inequality by Theorem~\ref{th:1}.

\begin{corollary}[(see {Patterson \cite{Pat}, Section 2.11})]
For any $t \in{\mathbb{R}}$ and $1/2 \le\sigma$, we have
\begin{eqnarray*}
\biggl| (\sigma+ \mathrm{{i}}t) (\sigma-1 + \mathrm{{i}}t) \uppi ^{-(\sigma
+\mathrm{{i}}t)/2} \Gamma
\biggl( \frac{\sigma+\mathrm{{i}}t}{2} \biggr) \zeta(\sigma +\mathrm{{i}}t) \biggr| \le \sigma(
\sigma-1) \uppi^{-\sigma/2} \Gamma \biggl( \frac{\sigma}{2} \biggr) \zeta(
\sigma). \label{eq:cor1}
\end{eqnarray*}
\end{corollary}

\subsection{Proof of Theorem \texorpdfstring{\protect\ref{th:2}}{1.2}}
Recall that ${\mathcal{Z}}$ is the set of zeros of the Riemann zeta
function which lie in the critical strip $\{s\in{\mathbb{C}} \dvt 0 <
\Re
(s) <1\}$ (see Section~\ref{sec1.3}). Observe that by the functional equation
and $\overline{\zeta(s)} = \zeta(\overline{s})$ if $\rho\in
{\mathcal
{Z}}$ then $\overline{\rho}, 1-\rho, 1-\overline{\rho} \in
{\mathcal
{Z}}$. There are no real elements of ${\mathcal{Z}}$ since $\xi
(\sigma
) < 0$ and $0< \Gamma(\sigma/2)$ when $0 < \sigma<1$ (see Section~\ref{sec1.1}
and the proof of Theorem~\ref{th:1}). Now we quote the following fact
from Patterson \cite{Pat}.

\begin{lemma}[(see {Patterson \cite{Pat}, page 34})]\label{lem:had}
Let ${\mathcal{Z}}_+ := \{ \rho\in{\mathcal{Z}} \dvt  \Im(\rho) >0 \}$.
Then $\sum_{\rho\in{\mathcal{Z}}_+} |\rho|^{-a}$ converges for all
$a>1$ and it holds that
%
\begin{equation}
\label{eq:Had} \xi(s) = s(s-1)\uppi^{-s/2} \Gamma \biggl( \frac{s}{2}
\biggr) \zeta (s) = \prod_{\rho\in{\mathcal{Z}}_+} \biggl( 1 -
\frac{s}{\rho} \biggr) \biggl( 1 - \frac{s}{1-\rho} \biggr)
\end{equation}
the product being absolutely convergent for all $s \in{\mathbb{C}}$.
\end{lemma}

\begin{pf*}{Proof of Theorem~\ref{th:2}}
If $\Xi_\sigma(t)$ is a pretended-infinitely divisible characteristic
function for any $1/2 < \sigma<1$, then $\zeta(s) \ne0$ for any $1/2
< \sigma<1$ by $\exp(z) \ne0$ for all $z \in{\mathbb{C}}$, $\Gamma
(s) \ne0$ for any $1/2 < \sigma<1$ and the representation (\ref{INF}).

Next suppose that the Riemann hypothesis is true. Then we have $\rho=
1/2 + \mathrm{{i}}\gamma$ and $1-\rho= 1/2 - \mathrm{{i}}\gamma$,
where $\gamma>0$ for $\rho\in{\mathcal{Z}}_+$. Note that the
exponential distribution with parameter $a>0$ is defined by $\mu(B) :=
a \int_{B \cap(0,\infty)} \mathrm{e}^{-ax}\,\mathrm{d}x$, where $B \in{\mathfrak{B}}
({\mathbb{R}})$. The characteristic function is given by $\widehat
{\mu
}(t)= a/(a-\mathrm{{i}}t)$ (see, for example, Sato \cite
{S99}, page 13).
Moreover, it is well known that
%
\begin{equation}
\label{eq:sato} \frac{a}{a-\mathrm{{i}}z} = \exp \biggl[ \int_0^\infty
\bigl(\mathrm{e}^{\mathrm
{i}zx}-1 \bigr) x^{-1} \mathrm{e}^{-a x} \,\mathrm{d}x \biggr],\qquad  a
>0 , z \in{\mathbb{R}}
\end{equation}
(see, for instance, Sato \cite{S99}, page 45). The formula
above holds if $a$
is replaced by $\alpha$ with $\Re(\alpha) >0$. This is proved as
follows. Put $\alpha= a+\mathrm{{i}}b$, $a>0$ and $b \in{\mathbb{R}}$.
Then one has
%
\begin{eqnarray}
\label{eq:gexp} \frac{\alpha}{\alpha-\mathrm{{i}}z}& =& \frac{a+\mathrm{{i}}b}{a} \frac
{a}{a+\mathrm{{i}}b-\mathrm{{i}}z}\nonumber
\\
&=& \exp \biggl[ \int_0^\infty
\bigl(\mathrm{e}^{\mathrm{i}(z-b)x}-1 \bigr) x^{-1} \mathrm{e}^{-a x} \,\mathrm{d}x -\int
_0^\infty \bigl( \mathrm{e}^{-\mathrm{i}bx}-1 \bigr)
x^{-1} \mathrm{e}^{-a x} \,\mathrm{d}x \biggr]
\\
&=& \exp \biggl[ \int_0^\infty
\bigl(\mathrm{e}^{\mathrm{i}zx}-1 \bigr) x^{-1} \mathrm{e}^{-\alpha x} \,\mathrm{d}x \biggr],\qquad  \Re(
\alpha) >0,\nonumber
\end{eqnarray}
by (\ref{eq:sato}). Thus, it holds that
\begin{eqnarray*}
\biggl( 1 - \frac{\sigma- \mathrm{{i}}t}{\rho} \biggr) \biggl( 1 - \frac{\sigma}{\rho}
\biggr)^{-1} &=& \frac{1/2-\sigma+ \mathrm{{i}}(\gamma+t)}{1/2+\mathrm{{i}}\gamma} \frac{1/2+\mathrm{{i}}\gamma}{1/2-\sigma+ \mathrm{{i}}\gamma} =
\frac{\sigma- 1/2 - \mathrm{{i}}\gamma-\mathrm{{i}}t}{\sigma- 1/2 -
\mathrm{{i}}\gamma}
\\
&=& \exp \biggl[ - \int_0^\infty
\bigl(\mathrm{e}^{\mathrm{i}tx}-1 \bigr) \mathrm{e}^{(1/2 -\sigma+
\mathrm{{i}}\gamma)x} \frac{\mathrm{d}x}{x} \biggr],
\end{eqnarray*}
where $\sigma>1/2$. It should be noted that we have $\sigma-\mathrm
{{i}}t \ne\rho, 1-\rho$ when $\sigma>1/2$ under the Riemann
hypothesis. Therefore, one has
%
\begin{eqnarray}
\label{eq:var} \varphi_\rho(t) &:=& \biggl( 1 - \frac{\sigma-\mathrm{{i}}t}{\rho}
\biggr) \biggl( 1 - \frac
{\sigma}{\rho} \biggr)^{-1} \biggl( 1 -
\frac{\sigma-\mathrm{{i}}t}{1-\rho} \biggr) \biggl( 1 - \frac{\sigma}{1-\rho} \biggr)^{-1}\nonumber
\\
&=& \frac{\sigma- 1/2 - \mathrm{{i}}\gamma-\mathrm{{i}}t}{\sigma-1/2 -
\mathrm{{i}}\gamma} \frac{\sigma- 1/2 + \mathrm{{i}}\gamma-\mathrm{{i}}t}{\sigma-1/2 +
\mathrm{{i}}\gamma} \\
&=& \exp \biggl[-2 \int
_0^{\infty} \bigl(\mathrm{e}^{\mathrm{i}tx}-1 \bigr)
\frac{\cos
(\gamma
x)}{x \mathrm{e}^{(\sigma-1/2)x}} \,\mathrm{d}x \biggr].\nonumber
\end{eqnarray}
We remark that $x^{-1} \cos(\gamma x) \mathrm{e}^{(1/2-\sigma)x} (\mathrm{d}x)$ is not a
measure but a signed measure since one has $-1 \le\cos(\gamma x) \le
1$ when $\gamma\in{\mathbb{R}}$. By (\ref{eq:Had}) and the definition
of $\Xi_\sigma(t)$, we have
\begin{eqnarray*}
\Xi_\sigma(t) &=& \prod_{\gamma\in{\mathcal{Z}}_+}
\frac{\sigma-1/2 - \mathrm
{{i}}\gamma+\mathrm{{i}}t}{\sigma-1/2 - \mathrm{{i}}\gamma} \frac{\sigma-1/2 + \mathrm{{i}}\gamma+\mathrm{{i}}t}{\sigma-1/2 +
\mathrm{{i}}\gamma}
\\
& =& \exp \biggl[ - 2 \sum_{1/2+\mathrm{{i}}\gamma\in{\mathcal{Z}}_+} \int
_0^{\infty} \bigl(\mathrm{e}^{\mathrm{i}tx}-1 \bigr)
\frac{\cos(\gamma x)}{x \mathrm{e}^{(\sigma-1/2)x}}\, \mathrm{d}x \biggr] .
\end{eqnarray*}
This equality implies (\ref{eq:th2}).
\end{pf*}

\begin{remark}
It should be mentioned that $\varphi_{1/2+\mathrm{{i}}\gamma} (t)$
defined by (\ref{eq:var}) is not a characteristic function for any
$\sigma>1/2$. It is proved by as follows. Obviously, one has
\[
\bigl |\varphi_{1/2+\mathrm{{i}}\gamma} (t)\bigr  |^2 = \frac{(\sigma-1/2)^2+ \gamma^2-t^2+(2\sigma-1) \mathrm
{{i}}t}{(\sigma
-1/2)^2 + \gamma^2} .
\]
If we take $t^2=2((\sigma-1/2)^2+ \gamma^2)$, then $|\varphi
_{1/2+\mathrm{{i}}\gamma} (t) |^2>1$.
\end{remark}

\subsection{Proof of Theorem \texorpdfstring{\protect\ref{th:3}}{1.3}}
Recall that ${\mathcal{Z}}$, ${\mathcal{Z}}_+^R$ and ${\mathcal
{Z}}_+^N$ is the set of zeros of $\zeta(s)$ which lie in $\{s\in
{\mathbb{C}} \dvt  0 < \Re(s) <1\}$, $\{s\in{\mathbb{C}} \dvt \Re(s) =1/2,
\Im(s) >0 \}$ and $\{s\in{\mathbb{C}} \dvt  1/2 < \Re(s) <1, \Im(s) >0
\}$, respectively. Then one has ${\mathcal{Z}}= \{ \rho, 1-\rho\dvt
\rho
\in{\mathcal{Z}}_+^R \} \cup\{ \rho, 1-\rho, \overline{\rho},
1-\overline{\rho}\dvt  \rho\in{\mathcal{Z}}_+^N \}$. We have the
following by Lemma~\ref{lem:had}.

\begin{lemma}
The sums $\sum_{\rho\in{\mathcal{Z}}_+^R} |\rho|^{-a}$ and $\sum_{\rho\in{\mathcal{Z}}_+^N} |\rho|^{-a}$ converge for all $a>1$ and
it holds that
%
\begin{eqnarray}
\label{eq:Had2} \xi(s) &= &\prod_{1/2+\mathrm{{i}}\gamma\in{\mathcal{Z}}_+^R} \biggl( 1 -
\frac{s}{1/2+\mathrm{{i}}\gamma} \biggr) \biggl( 1 - \frac
{s}{1/2-\mathrm{{i}}\gamma} \biggr)\nonumber
\\[-8pt]\\[-8pt]
&&{} \times \prod_{\rho\in{\mathcal{Z}}_+^N} \biggl( 1 - \frac{s}{\rho}
\biggr) \biggl( 1 - \frac{s}{1-\rho} \biggr) \biggl( 1 - \frac{s}{\overline{\rho}}
\biggr) \biggl( 1 - \frac
{s}{1-{\overline{\rho}}} \biggr),\nonumber
\end{eqnarray}
the products being absolutely convergent for all $s \in{\mathbb{C}}$.
\end{lemma}

\begin{pf*}{Proof of Theorem~\ref{th:3}}
Put $\overline{s} = \sigma-\mathrm{{i}}t$. Then we have
\begin{eqnarray*}
&& \biggl( 1 - \frac{\overline{s}}{\rho} \biggr) \biggl( 1 - \frac
{\sigma
}{\rho}
\biggr)^{-1} \biggl( 1 - \frac{\overline{s}}{1 - \rho} \biggr) \biggl( 1 -
\frac
{\sigma}{1 - \rho} \biggr)^{-1} \\
&&\quad = \frac{\sigma- \beta- \mathrm{{i}}\gamma-\mathrm{{i}}t}{\sigma-
\beta- \mathrm{{i}}\gamma} \frac{\sigma- 1 + \beta- \mathrm{{i}}\gamma-\mathrm{{i}}t}{\sigma
- 1+\beta- \mathrm{{i}}\gamma}
\\
&&\quad = \exp \biggl[ -\int_0^\infty
\bigl(\mathrm{e}^{\mathrm{i}tx}-1 \bigr) \mathrm{e}^{(\beta
-\sigma+ \mathrm{{i}}\gamma)x} \frac{\mathrm{d}x}{x} - \int
_0^\infty \bigl(\mathrm{e}^{\mathrm{i}tx}-1 \bigr)
\mathrm{e}^{(1-\beta-\sigma+
\mathrm{{i}}\gamma)x} \frac{\mathrm{d}x}{x} \biggr]
\end{eqnarray*}
from (\ref{eq:gexp}). By replacing $\rho$ by $\overline{\rho}$, we obtain
\begin{eqnarray*}
&& \biggl( 1 - \frac{\overline{s}}{\overline{\rho}} \biggr) \biggl( 1 - \frac{\sigma}{\overline{\rho}}
\biggr)^{-1} \biggl( 1 - \frac{\overline{s}}{1-\overline{\rho}} \biggr) \biggl( 1 -
\frac{\sigma}{1-\overline{\rho}} \biggr)^{-1}
\\
&&\quad = \exp \biggl[ - \int_0^\infty
\bigl(\mathrm{e}^{\mathrm{i}tx}-1 \bigr) \mathrm{e}^{(\beta-\sigma-
\mathrm{{i}}\gamma)x} \frac{\mathrm{d}x}{x} - \int
_0^\infty \bigl(\mathrm{e}^{\mathrm{i}tx}-1 \bigr)
\mathrm{e}^{(1-\beta-\sigma-
\mathrm{{i}}\gamma)x} \frac{\mathrm{d}x}{x} \biggr].
\end{eqnarray*}
We have to mention that one has $\beta-\sigma<0$ and $1-\beta-\sigma
<0$ since $\zeta(s) \ne0$ for $\sigma\ge1$ (see Remark~\ref
{re:zerofree} below). Hence, one has
\begin{eqnarray*}
&&\frac{(1-\overline{s}/\rho)(1-\overline{s}/\overline{\rho})
(1-\overline{s}/(1 - \rho) ) (1-\overline{s}/(1 -
\overline
{\rho}) )} {
(1-\sigma/\rho)(1-\sigma/\overline{\rho}) (1-\sigma/(1 -
\rho)
) (1-\sigma/(1 - \overline{\rho}) )}
\\
&&\quad =  \exp \biggl[ - 2 \int_0^{\infty}
\bigl(\mathrm{e}^{\mathrm
{i}tx}-1 \bigr) \cos(\gamma x) \bigl(\mathrm{e}^{(\beta-\sigma)x} +
\mathrm{e}^{(1-\beta-\sigma)x} \bigr) \frac{\mathrm{d}x}{x} \biggr] .
\end{eqnarray*}
Therefore, we have
\begin{eqnarray*}
\Xi_\sigma(t) &=& \exp \biggl[  - 2 \sum_{1/2+\mathrm{{i}}\gamma\in{\mathcal{Z}}_+^R}
\int_0^{\infty} \bigl(\mathrm{e}^{\mathrm{i}tx}-1 \bigr) \cos(
\gamma x) \mathrm{e}^{(1/2-\sigma)x} \frac{\mathrm{d}x}{x}
\\
&&\hphantom{\exp \biggl[}{}- 2 \sum_{\beta
+\mathrm{{i}}\gamma\in{\mathcal{Z}}_+^N} \int_0^{\infty}
\bigl(\mathrm{e}^{\mathrm{i}tx}-1 \bigr) \cos(\gamma x) \bigl(\mathrm{e}^{(\beta-\sigma)x} +
\mathrm{e}^{(1-\beta-\sigma)x} \bigr)\frac{\mathrm{d}x}{x} \biggr]
\end{eqnarray*}
by (\ref{eq:Had2}) and the definition of $\Xi_\sigma(t)$.
\end{pf*}

\begin{remark}\label{re:zerofree}
By modifying the proof above, we can see that one has (\ref{eq:th3})
for any $\sigma\ge\sigma_0 > 1/2$ if $\zeta(s)$ does not vanish for
$\sigma\ge\sigma_0$.
\end{remark}

\subsection{Proof of Theorem \texorpdfstring{\protect\ref{th:4}}{1.4}}
In order to prove Theorem~\ref{th:4}, we first prove the following
lemma which is an analogue of Nikeghbali and Yor \cite{NiYo}, Lemma~2.9.

\begin{lemma}\label{lem:ga}
Let $G_\sigma(t) = \Gamma(\sigma-\mathrm{{i}}t) / \Gamma(\sigma)$
for $0< \sigma$. Then $G_\sigma(t)$ is an infinitely divisible
characteristic function for any $\sigma>0$. Moreover, one has
\begin{eqnarray*}
\log G_\sigma(t) &=& \mathrm{i}t \lambda_\sigma^{\#}
+ \int_0^\infty \bigl( \mathrm{e}^{\mathrm{{i}}tx}-1-
\mathrm{{i}}tx 1_{[0,1]}(x) \bigr)\nu_{\sigma}^{\#}(\mathrm{d}x),
\\
\lambda_\sigma^{\#} &=& C(\sigma) := \int_0^1
\biggl(\frac
{\mathrm{e}^{-\sigma
x}}{1-\mathrm{e}^{-x}} - \frac{\mathrm{e}^{-x}}{x} \biggr) \,\mathrm{d}x - \int
_1^\infty \mathrm{e}^{-x} \frac{\mathrm{d}x}{x},
\\
\nu_{\sigma}^{\#}(\mathrm{d}x)&:=& \frac{1(\mathrm{d}x)}{x\mathrm{e}^{\sigma x}(1-\mathrm{e}^{-x})} .
\end{eqnarray*}
\end{lemma}

\begin{pf}
By the integral representation of $\Gamma(s)$ and the change of
variables $x=\mathrm{e}^{-y}$, we have
\begin{eqnarray*}
G_\sigma(t) &=& \frac{1}{\Gamma(\sigma)} \int_0^{\infty}
\mathrm{e}^{-x}x^{\sigma-1-\mathrm{{i}}t}\,\mathrm{d}x = \frac{-1}{\Gamma(\sigma)} \int_{\infty}^{-\infty}
\mathrm{e}^{-\mathrm{e}^{-y}} \mathrm{e}^{y(1-\sigma+\mathrm{{i}}t)y} \mathrm{e}^{-y}\,\mathrm{d}y
\\
&= & \frac{1}{\Gamma(\sigma)} \int_{-\infty}^{\infty}
\mathrm{e}^{\mathrm
{{i}}ty} \exp \bigl( -\sigma y -\mathrm{e}^{-y} \bigr) \,\mathrm{d}y,\qquad  \sigma>0.
\end{eqnarray*}
Therefore, the probability density function is given by $\exp( -\sigma
y -\mathrm{e}^{-y} )/\Gamma(\sigma)$.

Next, we quote Malmst\'en's formula (see, for example, Whittaker and Watson
\cite{WW}, page 249)
\[
\log\Gamma(s) = \int_0^\infty \biggl(
\frac{\mathrm{e}^{-sx}-\mathrm{e}^{-x}}{1-\mathrm{e}^{-x}} + (s-1) \mathrm{e}^{-x} \biggr) \frac{\mathrm{d}x}{x}, \qquad \sigma>0.
\]
Hence, it holds that
\begin{eqnarray*}
\log G_\sigma(t) &=& \int_0^\infty \biggl(
\frac{\mathrm{e}^{-(\sigma-\mathrm{{i}}t)x}-\mathrm{e}^{-\sigma x}}{1-\mathrm{e}^{-x}} - \mathrm{{i}}t \mathrm{e}^{-x} \biggr) \frac{\mathrm{d}x}{x}
\\
&=& \int_0^1 \biggl( \frac{\mathrm{e}^{\mathrm{{i}}tx} -1-\mathrm
{{i}}tx}{\mathrm{e}^{\sigma
x}(1-\mathrm{e}^{-x})} -
\mathrm{{i}}t \mathrm{e}^{-x} + \frac{\mathrm{{i}}tx\mathrm{e}^{-\sigma x}}{1-\mathrm{e}^{-x}} \biggr) \frac{\mathrm{d}x}{x} +
\int_1^\infty \biggl( \frac{\mathrm{e}^{\mathrm{{i}}tx} -1}{\mathrm{e}^{\sigma x}(1-\mathrm{e}^{-x})} -
\mathrm{{i}}t \mathrm{e}^{-x} \biggr) \frac{\mathrm{d}x}{x}
\\
&=& \int_0^\infty\frac{\mathrm{e}^{\mathrm{{i}}tx}-1-\mathrm{{i}}tx1_{[0,1]}(x)}{x
\mathrm{e}^{\sigma x}(1-\mathrm{e}^{-x})}\,\mathrm{d}x +
\mathrm{{i}}t \int_0^1 \biggl(
\frac{\mathrm{e}^{-\sigma x}}{1-\mathrm{e}^{-x}} - \frac
{\mathrm{e}^{-x}}{x} \biggr)\, \mathrm{d}x - \mathrm{{i}}t \int
_1^\infty \mathrm{e}^{-x} \frac{\mathrm{d}x}{x} .
\end{eqnarray*}
Therefore, we obtain Lemma~\ref{lem:ga}.
\end{pf}

For the reader's convenience, we give a proof of (\ref{eq:rzlm1}). By
the Euler product of $\zeta(s)$ and the Taylor expansion of $\log
(1-x)$, $|x|<1$, one has
\begin{eqnarray*}
 \log\frac{\zeta(\sigma- \mathrm{i}t)}{\zeta(\sigma)} &=& \sum_p \log
\frac{1-p^{-\sigma}}{1-p^{-\sigma+\mathrm{i}t}}= \sum_p \sum
_{r=1}^\infty\frac{1}{r}p^{-r\sigma} \bigl(
p^{r\mathrm
{i}t}-1 \bigr)
\\
&=& \sum_p \sum_{r=1}^\infty
\frac{1}{r}p^{-r\sigma} \bigl( \mathrm{e}^{r\mathrm
{i}t\log p}-1 \bigr) = \int
_{-\infty}^{\infty} \bigl( \mathrm{e}^{\mathrm{i}tx}-1 \bigr) \sum
_{p}\sum_{r=1}^{\infty}
\frac{1}{r}p^{r\sigma}\delta_{r\log
p}(\mathrm{d}x) .
\end{eqnarray*}
This equality implies (\ref{eq:rzlm1}).
\begin{pf*}{Proof of Theorem~\ref{th:4}}
We have
\[
\Xi_\sigma(t) = \uppi^{\mathrm{i}t/2} G_{\sigma/2} (t/2)
\frac{\sigma- \mathrm{{i}}t}{\sigma} \frac{\sigma-1 - \mathrm
{{i}}t}{\sigma-1} \frac{\zeta(\sigma-\mathrm{i}t)}{\zeta(\sigma)}
\]
by the definition of $\Xi_\sigma(t)$. It holds that
\begin{eqnarray*}
\log G_{\sigma/2} (t/2) &=& \frac{\mathrm{{i}}t}{2} C(\sigma/2) + \int
_0^\infty\frac{\mathrm{e}^{\mathrm{{i}}(t/2)x}-1-\mathrm
{{i}}(t/2)x1_{[0,1]}(x)}{x \mathrm{e}^{\sigma x/2}(1-\mathrm{e}^{-x})}\,\mathrm{d}x
\\
&=& \mathrm{{i}}t \frac{C(\sigma/2)}{2} + \int_0^\infty
\frac{\mathrm{e}^{\mathrm{{i}}tx}-1-\mathrm
{{i}}tx1_{[0,1/2]}(x)}{x \mathrm{e}^{\sigma x}(1-\mathrm{e}^{-2x})}\,\mathrm{d}x
\end{eqnarray*}
from Lemma~\ref{lem:ga}. Obviously, one has $1/2<r\log p$ for any
integer $r$ and prime number $p$ since $\log2 = 0.6931471806\ldots\,$.
Hence by using (\ref{eq:rzlm1}), we have
\[
\log\frac{\zeta(\sigma-\mathrm{i}t)}{\zeta(\sigma)} = \int_{-\infty}^{\infty} \bigl(
\mathrm{e}^{\mathrm{i}tx}-1 -\mathrm {{i}}tx1_{[0,1/2]}(x) \bigr) \sum
_{p}\sum_{r=1}^{\infty}
\frac{1}{r}p^{-r\sigma}\delta_{r\log
p}(\mathrm{d}x) .
\]
When $\sigma>1$, one has
\begin{eqnarray*}
\frac{\sigma- \mathrm{{i}}t}{\sigma} &=&\exp \biggl[- \int_0^\infty
\bigl(\mathrm{e}^{\mathrm{i}tx}-1-\mathrm{{i}}tx1_{[0,1/2]}+\mathrm
{{i}}tx1_{[0,1/2]} \bigr) \mathrm{e}^{-\sigma x} \frac{\mathrm{d}x}{x} \biggr]
\\
&=& \exp \biggl[- \int_0^\infty\frac{\mathrm{e}^{\mathrm{i}tx}-1-\mathrm
{{i}}tx1_{[0,1/2]}(x)}{x \mathrm{e}^{\sigma x}}
\,\mathrm{d}x - \mathrm{i}t \frac{1-\mathrm{e}^{-\sigma/2}}{\sigma} \biggr]
\end{eqnarray*}
by (\ref{eq:sato}). Thus, it holds that
\begin{eqnarray*}
&&\frac{\sigma- \mathrm{{i}}t}{\sigma} \frac{\sigma-1 - \mathrm
{{i}}t}{\sigma-1}
\\
&&\quad = \exp \biggl[- \int_0^\infty \bigl(
\mathrm{e}^{\mathrm{i}tx}-1-\mathrm {{i}}tx1_{[0,1/2]}(x) \bigr) \frac{1+\mathrm{e}^x}{x \mathrm{e}^{\sigma x}}
\,\mathrm{d}x - \mathrm{i}t \biggl( \frac{1-\mathrm{e}^{-\sigma
/2}}{\sigma} + \frac{1-\mathrm{e}^{-(\sigma-1)/2}}{\sigma-1} \biggr) \biggr].
\end{eqnarray*}
If $x$ is sufficiently large, then we have
\[
\frac{1}{x \mathrm{e}^{\sigma x}(1-\mathrm{e}^{-2x})}-\frac{1+\mathrm{e}^x}{x\mathrm{e}^{\sigma x}} < 0 .
\]
Thus $\nu_\sigma$ in Theorem~\ref{th:4} is not a measure but a
signed measure.

Finally, we show $\int_{\mathbb{R}} (|x|^{2} \wedge1) |\nu_\sigma|(\mathrm{d}x)
< \infty$ when $\sigma>1$. By using $(1-\mathrm{e}^{-2})x \le1-\mathrm{e}^{-2x}$ for $0
\le x <1$ and $1-\mathrm{e}^{-2} \le1-\mathrm{e}^{-2x}$ for $x \ge1$, we have
\[
\int_0^\infty\frac{(1-\mathrm{e}^{-2})(|x|^{2} \wedge1)}{x \mathrm{e}^{\sigma
x}(1-\mathrm{e}^{-2x})} \,\mathrm{d}x \le \int
_0^1 \frac{\mathrm{d}x}{\mathrm{e}^{\sigma x}} + \int
_1^\infty\frac
{\mathrm{d}x}{x\mathrm{e}^{\sigma
x}} < \int
_0^\infty\frac{\mathrm{d}x}{\mathrm{e}^{\sigma x}} < \infty.
\]
Obviously, it holds that
\[
\int_0^\infty\frac{(1+\mathrm{e}^x)(|x|^{2} \wedge1)}{x\mathrm{e}^{\sigma x}}\, \mathrm{d}x < 2 \int
_0^\infty\frac{(|x|^{2} \wedge1)\,\mathrm{d}x}{x\mathrm{e}^{(\sigma-1)x}} < 2 \int
_0^\infty\frac{\mathrm{d}x}{\mathrm{e}^{(\sigma-1)x}} < \infty.
\]
From $\sum_{p} p^{-\sigma} < \sum_{n=2}^\infty n^{-\sigma} = \zeta
(\sigma) -1$, one has
\begin{eqnarray*}
&&\int_0^\infty\sum_{p}
\sum_{r=1}^{\infty}\frac{1}{r}p^{-r\sigma
}
\delta _{r\log p}(\mathrm{d}x) \\
&&\quad = \sum_{p}\sum
_{r=1}^{\infty}\frac{1}{r}p^{-r\sigma} < \sum
_{p}\sum_{r=1}^{\infty}p^{-r\sigma}< \sum_{n=1}^\infty n^{-\sigma} + \sum
_{p}\sum_{r=2}^{\infty
}p^{-r\sigma}
\\
&&\quad
= \zeta(\sigma) + \sum_{p} \frac{p^{-2\sigma}}{1-p^{-\sigma}} <
\zeta (\sigma) + \sum_{n=2}^\infty
\frac{n^{-2\sigma}}{1-2^{-\sigma}}
\\
&&\quad  < \zeta(\sigma) + \bigl(1-2^{-\sigma}\bigr)^{-1} \zeta(2\sigma)
< \infty.
\end{eqnarray*}
Therefore the characteristic function $\Xi_\sigma(t)$ is not
infinitely divisible but quasi-infinitely divisible.
\end{pf*}

\begin{remark}
Suppose $\sigma\ne1$ and put
\[
\Xi_\sigma^* (t) := \frac{\sigma-1}{\sigma-1 - \mathrm{{i}}t} \Xi _\sigma(t).
\]
Then $\Xi_\sigma^*(t)$ is a characteristic function for any $\sigma
\ne
1$ by the fact that the product of a finite number of characteristic
functions is also a characteristic function. By modifying the proof
above, we have
\begin{eqnarray*}
\Xi_\sigma^* (t)& =& \exp \biggl[ \mathrm{i}t \lambda_\sigma^*
+ \int_0^\infty \bigl(\mathrm{e}^{\mathrm{{i}}tx}-1-\mathrm
{{i}}tx1_{[0,1/2]}(x) \bigr)\nu_{\sigma}^*(\mathrm{d}x) \biggr],
\\
\lambda_\sigma^* &:=& \frac{1-\mathrm{e}^{-\sigma/2}}{\sigma} + \frac{\log
\uppi
}{2} +
\frac{1}{2} \int_0^1 \biggl(
\frac{\mathrm{e}^{-\sigma x/2}}{1-\mathrm{e}^{-x}} - \frac
{\mathrm{e}^{-x}}{x} \biggr)\, \mathrm{d}x - \frac{1}{2} \int
_1^\infty \mathrm{e}^{-x} \frac{\mathrm{d}x}{x},
\\
\nu_{\sigma}^*(\mathrm{d}x)&:=& \frac{1(\mathrm{d}x)}{x \mathrm{e}^{\sigma x}(1-\mathrm{e}^{-2x})}-\frac
{1(\mathrm{d}x)}{x\mathrm{e}^{\sigma x}} + \sum
_{p}\sum_{r=1}^{\infty}
\frac{p^{-r\sigma}}{r}\delta_{r\log p}(\mathrm{d}x)
\end{eqnarray*}
for $\sigma>1$. Therefore the characteristic function $\Xi_\sigma^*
(t)$ is infinitely divisible for any $\sigma>1$ since one has
\[
\frac{1}{x \mathrm{e}^{\sigma x}(1-\mathrm{e}^{-2x})}-\frac{1}{x\mathrm{e}^{\sigma x}} > 0 ,\qquad  x>0.
\]

Moreover, we can see that every characteristic function $\Xi_\sigma
^*(t)$ is a pretended-infinitely divisible characteristic function for
each $1/2 < \sigma<1$ if and only if the Riemann hypothesis is true by
an argument similar to that in the proof of Theorem~\ref{th:2}. In
addition, it holds that
\begin{eqnarray*}
\Xi_\sigma^* (t) &= &\exp \biggl[ \int_{-\infty}^\infty
\bigl(\mathrm{e}^{\mathrm{{i}}tx}-1 \bigr)\nu _{\sigma}^*(\mathrm{d}x) \biggr],
\\
\nu_{\sigma}^*(\mathrm{d}x)&:=& \frac{1_{(-\infty,0)}(\mathrm{d}x)}{-x\mathrm{e}^{(\sigma-1)x}} -\sum
_{1/2+\mathrm{{i}}\gamma\in{\mathcal{Z}}_+} \frac{2\cos
(\gamma
x)}{x\mathrm{e}^{(\sigma-1/2)x}} 1_{(0,\infty)}(\mathrm{d}x) ,
\end{eqnarray*}
for $1/2 < \sigma<1$, under the Riemann hypothesis. This is proved by
(\ref{eq:th2}) and
\begin{eqnarray*}
\frac{\sigma-1}{\sigma-1 - \mathrm{{i}}t}& =& \frac{1-\sigma
}{1-\sigma
+ \mathrm{{i}}t} = \exp \biggl[ \int
_0^\infty \bigl(\mathrm{e}^{-\mathrm{i}tx}-1 \bigr)
\mathrm{e}^{(\sigma
-1)x} \frac{\mathrm{d}x}{x} \biggr]
\\
&=& \exp \biggl[ \int_0^{-\infty} \frac{(\mathrm{e}^{\mathrm{i}tx}-1)\,\mathrm{d}x}{x
\mathrm{e}^{(\sigma
-1)x}}
\biggr] = \exp \biggl[ - \int_{-\infty}^0
\frac{(\mathrm{e}^{\mathrm{i}tx}-1)\,\mathrm{d}x}{x
\mathrm{e}^{(\sigma-1)x}} \biggr]
\end{eqnarray*}
when $1/2 < \sigma<1$.

It is well known that convolving a density with a normal density to
make distributions more well-behaved. In this case the exponential
distribution is the one that makes things nicer since when $\sigma>1$,
the complete Riemann zeta distribution defined by $\Xi_\sigma(t)$ and
the distribution defined by $\Xi_\sigma^*(t)$ are quasi-infinitely
divisible and infinitely divisible, respectively.
\end{remark}

\section*{Acknowledgements}
The author would like to thank for Professor Takahiro Aoyama for his
useful comments. Moreover, the author is grateful to the referee for
his or her helpful comments and suggestions on the manuscript.




\printhistory

\end{document}